\documentclass[12pt]{article} %H

\usepackage{amsmath,amsfonts,amssymb}

\let\emptyset\varnothing
\def\dim{\mathop{\mathrm{dim}\;\!}\nolimits}
\def\codim{\mathop{\mathrm{codim}\;\!}\nolimits}

\def\hit{\mathbb{H}} %{ \hbox{ I\hskip -2pt R} }
\def\ait{\mathbb{A}}    %{ \hbox{\bf {\it  A\hskip -10pt A}} }
\def\ppit{\mathbb{P}} %{ \hbox{\bf {\it I\hskip -2pt P}} }
\def\Ai{\overline{A}}
\def\Hi{H_{\infty}}
\def\A{A}

\def\Ch{\mathsf{Ch}}
\def\bnbc{\beta{\textup{\textbf{nbc}}}}

\def\PM{\mathsf{PM}}
\def\th{\widetilde{h_i}}
\def\tg{\widetilde{g_{i,j}}}
\def\oC{\overline{C}}
\def\l{\mathbf{l}}
\def\m{\mathbf{m}}
\def\Re{\mathsf{Re}\,}
\def\pf{\mathsf{Proof.}}
\def\dind{\Delta,\, F_{tr(\Delta)}=F}
\def\ftext#1{{\let\thefootnote\relax\footnotetext{\noindent #1}}}

\newcommand{\al}{\alpha}
\newcommand{\om}{\omega}

\newcommand{\qed}{\hfill~~\mbox{$\square$}}

\newtheorem{theorem}{Theorem}[section]
\newtheorem{proposition}[theorem]{Proposition}
\newtheorem{definition}[theorem]{Definition}
\newtheorem{lemma}[theorem]{Lemma}
\newtheorem{corollary}[theorem]{Corollary}

\begin{document}

\title{\Large \bf The Determinant of a Hypergeometric Period Matrix}
\author{ \normalsize Y. Markov$^{\star}\,$, V. Tarasov$^*\,$
and A. Varchenko$^{\diamond}$}
\date{}
\maketitle
\begin{center}
{\it
$^{\star,\,\diamond}$ Department of Mathematics, University of North
Carolina,\\
Chapel Hill, NC 27599 -- 3250, USA\\ \hfill\\
$^{\ast}$ Department of Mathematics, Faculty of Science, Osaka
University\\
Toyonaka, Osaka 560, Japan}
\end{center}

\medskip
\centerline{September, 1997}
\medskip
\ftext{\hspace{-0.6cm}$^{\star}$ {\sl E-mail\/{\rm:} markov@math.unc.edu}\\
$^*$ {\sl E-mail\/{\rm:} vt@math.sci.osaka-u.ac.jp\,, \,vt@pdmi.ras.ru}\\
\hphantom{$*$} On leave of absense from St.\,Petersburg Branch of
Steklov Mathematical Institute\\
$^{\diamond}$ {\sl E-mail\/{\rm:} av@math.unc.edu}\\
\hphantom{$*$} Supported in part by NSF grant DMS-9501290
}

\section{Introduction}

The Euler beta function is an alternating product of Euler gamma functions,
\begin{equation}\label{2}
B (\alpha, \beta)\,=\, {\Gamma (\alpha)\,\Gamma (\beta)
\over
\Gamma (\alpha + \beta)}
\end{equation}
where the
 Euler gamma  and beta functions are defined by
\begin{equation}\label{1}
\Gamma (\alpha)\,=\,\int_0^\infty\, t^{\alpha -1}\, e^{-t}\,dt,
\qquad
B (\alpha, \beta)\,=\,\int_0^1\, t^{\alpha -1}\,(1-t)^{\beta -1}\, \,dt.
\end{equation}
There is a generalization of formula (1) to the case of an arrangement
of hyperplanes in an affine space, see [V1, V2, DT].

{\bf Example.} Consider an arrangement of three points $z_1, z_2, z_3$ in a
line. The point $z_j$ is the zero of the function $f_j=t-z_j$. Set
$$
\Delta_1= [z_1,z_2],
\qquad
\Delta_1= [z_2,z_3],
$$
$$
U_\al\,=\,(t-z_1)^{\alpha_1}(t-z_2)^{\alpha_2}(t-z_3)^{\alpha_3},
$$
$$
\om_1\,=\,\al_1 U_\al dt/(t-z_1), \qquad
\om_2\,=\,\al_2 U_\al dt/(t-z_2),
$$
then
$$
\det\Bigl(\;\!\int_{\Delta_i}\om_j\Bigr)\,=\,
{\Gamma (\alpha_1+1)\,\Gamma (\alpha_2+1)\,\Gamma (\alpha_3+1)
\over
\Gamma (\alpha_1+\alpha_2+\alpha_3+1)}\,\prod_{i\neq j} f_i^{\al_i}(z_j).
$$

In this paper we describe a generalization of the first formula in (2)
to the case of an arrangement
of hyperplanes in an affine space.

{\bf Example.} Consider an arrangement of two points $z_1, z_2$ in a
line. Let $f_j=t-z_j$. Set
$$
\Delta_1= [z_1,z_2],
\qquad
\Delta_2= [z_2,\infty],
$$
$$
U_\al\,=\,(t-z_1)^{\alpha_1}(t-z_2)^{\alpha_2},
$$
$$
\om_1\,=\,\al_1 U_\al dt/(t-z_1), \qquad
\om_2\,=\,\al_2 U_\al dt/(t-z_2).
$$
Let $a$ be a positive number. Then
$$
\det\Bigl(\;\!\int_{\Delta_i}e^{-at}\om_j\Bigr)\,=\,
\Gamma (\alpha_1+1)\Gamma (\alpha_2+1) \, e^{-a(z_1 + z_2)}\,
a^{-(\alpha_1+\alpha_2)}\,\prod_{i\neq j} f_i^{\al_i}(z_j).
$$

The determinant formulas  are useful, in particular
in applications to the Knizhnik-Zamolodchikov type of differential equations
when a determinant formula allows one to conclude that a set of solutions
to the equation given by suitable multidimensional hypergeometric
integrals forms a basis of solutions, cf. [SV], [TV], [V3], see also [L],
[LS], [V4], [V5].
\vskip\baselineskip

The paper is organized as follows. Sections~\ref{arrangements} --
\ref{hg-section} contain definitions of the main objects: arrangements,
critical values, and hypergeometric period matrices. The main result of
the paper is Theorem~\ref{maintheorem}. The proofs of all statements are
presented in Section~\ref{proofs}. In Sections~\ref{selberg} and
\ref{pfselb} we discuss two determinant formulas for Selberg type
integrals. In this case the configuration of hyperplanes is special and
admits a symmetry group. The symmetry group acts on the domains of the configuration
and on the hypergeometric differential forms associated with the configuration.
Therefore the period matrix of the configuration $( \,\int_{\Delta_i} \omega_j\,)$
splits into blocks according to different representations of the symmetry group.
We compute the determinant of the block corresponding to
the trivial representation.

\section{Arrangements}
\label{arrangements}

In this section we review results from [FT] and [V1].

\subsection{}
 Let $f_{1},\ldots ,f_{p}$ be linear polynomials  on a
real affine space $V$.
 Let $I$ denote $\{1,\ldots ,p\}$ and let $A$ be the arrangement
$\{H_{i}\}_{i\in I}$, where $H_{i}= \ker f_{i}$ is the hyperplane defined
 by $f_{i}$.

An {\it edge} of $A$ is a nonempty intersection of some of its
hyperplanes. A {\it vertex} is a \mbox{$0$-}dimensional edge.
Let $L(A)$ denote the set of all  edges.

An arrangement  $A$ is said to be {\it essential} if it has vertices.
Until the end of this paper we suppose that $A$ is essential.

An arrangement $A$ is said to be in {\it general position} if,
for all  subarrangements $\{H_{i_{1}},\cdots ,H_{i_{k}}\}$ of $A$, we have
$\codim(H_{i_{1}}\cap\cdots\cap H_{i_{k}})=k$ if $1\leq k\leq\dim V$ and
$H_{i_{1}}\cap\cdots\cap H_{i_{k}}=\emptyset$ if $k> \dim V$.

 Let
\begin{equation} \label{defM}
M(A)=V-\cup_{i\in I} H_{i}.
\end{equation}

The topological space $M (A)$ has finitely many connected components,
which are called {\it domains}. Domains are open   polyhedra, not
necessary bounded. Their faces are  precisely the domains of the
arrangements induced by $A$ on the edges of  $A$. More generally, in any
subspace $U \subset V$ the arrangement  $A$  cuts out a new arrangement
$A_U$ consisting of  the hyperplanes $\{ H_i\cap U \, | \, H_i \in A,
\, U \not \subset H_i\}$. $A_U$ is called a {\it section} of the
arrangement. For every edge $F$ of $A$ the domains of
the section $A_F$ are called the {\it faces} of the arrangement $A$.

Let $F$ be an edge of $A$ and $I(F)$ the set of all indices $i$ for
which $F \in H_i$. The arrangement $A^F$  in $V$ consisting of the
hyperplanes  $\{H_i \, | \, H_i \in A, \, i\in I(F)\}$, is called the
{\it localization of the arrangement at the edge $F$}.

Every edge $F$ of codimension $l$ is associated to an arrangement in
an \mbox{$(l-1)$--}dimensional projective space. Namely, let $L$ be a normal
subspace to $F$ of the complementary dimension.
 Consider the localization at this edge and its section by the
normal subspace. All of the hyperplanes of the resulting arrangement $(A^F)_L$
pass through the point $v=F\cap L$. We consider the arrangement which
$(A^F)_L$ induces in the tangent space $T_vL$. It determines an
arrangement in the projectivization of the tangent space, which is
called the {\it projective normal arrangement} and denoted $PA^F$.
The arrangements corresponding to different normal subspaces are naturally
isomorphic.

 A face of an arrangement is said to be {\it bounded relative to a
hyperplane} if the closure of the face does not intersect the hyperplane.
It is known [V1, Theorem 1.5] that if $A=\{H_i\}_{i\in I}$ is an arrangement
in a real projective space, then the number of domains bounded with respect
to $H_i$ does not depend on $i$. This number is called the {\it
discrete length of the arrangement}. The discrete length of the empty
arrangement is set to be equal to 1.

Let $F$ be an edge of an arrangement $A$ in a projective space.  The {\it
discrete length} of the edge is defined as the discrete length of the
arrangement $A_F$; the {\it discrete width} of the edge is the discrete
length of the arrangement $PA^F$; and the {\it discrete volume} of the
edge is the product of its discrete length and discrete width. These
numbers are denoted  $l(F)$, $s(F)$, and $vol(F)$, respectively.

If $F$ is a \mbox{$k$--}dimensional  edge of an arrangement $A$ in an affine
space, then its {\it
discrete length} is the number of bounded \mbox{$k$--}dimensional faces of the
arrangement $A_F$. Its {\it discrete width} is the discrete
length of the arrangement $PA^F$, and its {\it discrete volume}  is the
product of its discrete length and discrete width.

  Another more invariant definition of the above quantities could be
given as follows (see~[OT]).
  Consider the complexification and then the projectivization of the
affine space $V$. Denote it  $\ppit V$. Let $\hit_i=\{f_i=0\}_{i\in I}$
be  hyperplanes in $\ppit V$,$\quad\hit_{\infty}$ the infinite hyperplane
and $\ait$ the arrangement in $\ppit V$ defined by all these hyperplanes.
Let $\chi (\ait)$ denote the Euler characteristic of $\ppit V -
\cup_{i\in\overline{I}}\hit_i$, where $\overline{I} = 1,\ldots,p,\infty$.
If $F$ is an edge of $\ait$, then
$$l(F)=|\chi(\ait_F)|,\quad s(F)=|\chi(P\ait^F)|,\quad vol(F)=l(F)s(F).$$

\subsection{The beta-function of an arrangement}

An arrangement is called {\it weighted} if a complex number is assigned
to every hyperplane of the arrangement. The complex numbers are called
{\it weights}. The weight of a hyperplane $H_i$ is denoted
$\alpha_i$.

The {\it weight} of an edge $F$ of a weighted arrangement is
the sum, $\alpha (F)$,
of the weights of the hyperplanes which contain $F$.

 Let $A$ be a weighted arrangement in an  affine space $V$. Make
$V$ into a projective space by adding the hyperplane $H_{\infty}$ at
infinity: $ \overline{V} = V \cup H_{\infty}$. For all $i\in I$ denote
$\overline{H_i}$ the projective closure of $H_i$ in $\overline{V}$.
  Let $\Ai=\{\overline{H_i}\}_{i\in I}\cup\{\Hi\}$  be the corresponding
projective arrangement in  $\overline{V}$. The arrangement $\Ai$ is
called the {\it projectivization} of  $A$. Set $\alpha_{\infty} =
%%!!
-(\alpha_1+\cdots+\alpha_p)$.  Let $L_{-}$ denote the set of all edges at
infinity of the arrangement  $\Ai$ and $L_{+}$ the set consisting
of all the other edges.

\begin{definition}
$(i)$ Let the weights $\alpha_1,\ldots ,\alpha_p$ of the hyperplanes be
complex numbers with positive real part. The beta--function  of  an
affine arrangement $A$ is defined by
$$B(A;\alpha) = \prod_{F\in L_{+}} \Gamma (\alpha (F) +1)^{vol(F)}
\big/   \prod_{F\in L_{-}} \Gamma (-\alpha (F) +1)^{vol(F)} ; $$
$(ii)$ In addition, let $H_0$ be a hyperplane in $V$ and $\overline{H_0}$
its closure in $\overline{V}$. The  beta-function of an affine
arrangement $A$ relative to the hyperplane $H_0$ is  defined by
$$B(A;\alpha ;H_0) = \prod_{F\in L_{+}} \Gamma (\alpha (F) +1)^{vol(F)}
\big/ \prod_{F\in L_{-},\, F\subset \overline{H_0}} \Gamma (-\alpha (F)
+1)^{vol(F)}.$$
\end{definition}
{\sc Example.} Let $n=\dim V$.  For an arrangement $A$ of $p$
hyperplanes in general position the above formulas take the form
$$B(A;\alpha) = \left(\frac{\Gamma(\alpha_1 +1)\ldots\Gamma (\alpha_p +1)}
      {\Gamma (\alpha_1+\cdots+\alpha_p +1)}\right)^{\binom{p-2}{n-1}};
\qquad B(A;\alpha;f_0) =
\left(\Gamma(\alpha_1+1)\ldots\Gamma(\alpha_p+1)\right)^{\binom{p-2}{n-1}}.$$

\subsection{Trace of an arrangement at infinity relative to a hyperplane}
\label{traceinf}

 Assume that an additional non-constant linear function $f_0$ on $V$ is
given. Denote the hyperplane $\{f_0=0\}$ by $H_0$.

Let A be an affine arrangement in the affine space $V$ .
Consider the projectivized arrangement $\Ai$ in $\overline{V}$ and
its section $\Ai_{\Hi}$. The intersection of $\Ai$ with the affine
space
$W=\Hi - \overline{H_0} \cap \Hi$ is called the {\it trace of the
arrangement $A$ at
infinity} relative to the hyperplane $H_0$ and is denoted $tr(A)_{H_0}$.

\begin{lemma}
  If the affine arrangement $A$ is given by the linear functions
$\{f_i\}_{i\in I}$ on $V$, then the affine arrangement $tr(A)_{H_0}$ is
given by  the linear functions $\{h_i=f_i^0/f_0^0\, | \, i\in I;\,\,
f_i^0/f_0^0  \not = const \}$ on $W$,  where $f_i^0$ denotes the
homogeneous part of $f_i$. \qed
\end{lemma}

\section{Properties of an arrangement}
\label{properties}
In this section we examine the properties of the unbounded domains of an
arrangement $A$  on which an additional linear function $f_0$ tends to
$+\infty$.
\vskip\baselineskip

 Let $A$ be an arrangement in an affine space $V$.
Consider its projectivization $\Ai$ in the projective space $\overline{V}$.
Let $\Delta$ be an unbounded face of the arrangement
$A$. Take  the closure $\overline{\Delta}$ of $\Delta$ in $\overline{V}$.
Consider the intersection
$\overline{\Delta}\cap H_{\infty}$. It is a union of faces
of the arrangement $\Ai_{\Hi}$. There is a unique one of highest
dimension. We call it the {\it trace} of $\Delta$ at infinity and
denote $tr(\Delta)$.

An unbounded face $\Delta$ of $A$ is called a {\it growing face
with respect to } $f_0$ if $f_0(x)$ tends to $+\infty$
whenever $x$ tends to infinity in $\Delta$.
 A face at
infinity  $\Sigma$ is called a {\it bounded face at infinity with respect to }
$f_0$ if
$\overline{\Sigma}\cap \overline{H_0}$ is empty.
In other words, $\Sigma$ is a bounded face at infinity  if and only if it is
a bounded face of the affine arrangement $tr(A)_{H_0}$.

  If $\Sigma$ is a face of $A$, denote $F_{\Sigma}$ the unique edge of
the smallest  dimension which contains $\Sigma$. Define the discrete length,
the discrete width and the volume of the face as the same quantities for
the corresponding edge.

\begin{theorem} \label{bijection}
The  trace map from the unbounded faces  of  $A$ to the faces of $\Ai$ in
$H_{\infty}$ has the following properties:

$(i)$ The trace of a growing face is a bounded
face at infinity.

$(ii)$ For any bounded face at infinity, $\Sigma$, there exist exactly
$s(\Sigma )$ growing domains with trace $\Sigma$, where $s(\Sigma )$
denotes the discrete width of this face.

$(iii)$ The number of growing domains of $A$ is equal to the sum of the
volumes of all edges of $\Ai$ at infinity which do not lie in
$\overline{H_0}$:
$$ \mathrm{\#}\,\mathrm{ growing}\,\mathrm{ domains} = \sum_{F\subseteq
H_{\infty},\, F \not\subseteq \overline{H_0}} vol(F).$$
\end{theorem}
Theorem~\ref{bijection} is proved in Section~\ref{biject}

\section{Critical values}
\label{critical}

The aim of this section is to define the critical values of the functions
$f_1^{\alpha_1}, \ldots , f_p^{\alpha_p}$ on the bounded domains of an
arrangement $A$ and the critical values of the same functions, with
respect to an additional linear function $f_0$, on the bounded and growing
domains of $A$.

\medskip

Let an arrangement $A$ be given by linear functions
$\{f_i\}_{i\in I}$ and let $\alpha = \{\alpha_i\}_{i\in I}$ be a
corresponding set of weights.

For every $i\in I$, a face of the arrangement $A$ on which $f_i$ is
constant is called a {\it critical face} with respect to $f_i$ and the
value of $f_i$  on that face is called a {\it critical value}.
In particular, each vertex is a critical face for every function $f_i$.

Assume that a function $|f_i|$ is bounded on a face $\Sigma$ of $A$.
The subset of $\overline{\Sigma}$ on which $|f_i|$  attains its maximum
is a union of critical faces. Among them, there is a  unique one of
highest dimension. It is called the {\it external support}  of the face
$\Sigma$ with respect to $f_i$.

Denote $\Ch(A)$ the set of all bounded domains of $A$. Let $\beta(A) =
\left|\Ch(A)\right|$. Enumerate the bounded
domains by numbers $1,\ldots,\beta(A)$.
For every $i\in I$ and $j\in \{1,\ldots ,\beta (A)\}$, choose a branch
of the multi--valued function $f_i^{\alpha_i}$ on the domain $\Delta_{j}$
and denote it $g_{i,j}$. Let $\Sigma_{i,j}$ be the external support of
$\Delta_j$ with respect to $f_i$. Define the {\it extremal critical
value of the chosen branch $g_{i,j}$ on $\Delta_j$}  as the number
$c(g_{i,j},\Delta_j)=g_{i,j}(\Sigma_{i,j})$. Denote  $c(A;\alpha)$
the product of all extremal critical values of the chosen branches,
$$ c(A;\alpha)=\prod_{j=1}^{\beta(A)}\prod_{i\in I}c(g_{i,j},\Delta_j).$$

Assume that an additional non-constant linear function $f_0$ on $V$ is
given.  Let $\Delta$ be bounded or growing domain of $A$. Then $f_0$ is
bounded  from below on $\Delta$.  The subset of $\overline{\Delta}$ on
which $e^{-f_0}$  attains its maximum coincides with the subset of
$\overline{\Delta}$  where $f_0$ attains its minimum. This subset has a
unique face of  highest dimension; it is called the {\it support face} of
$f_0$ on $\Delta$ and denoted $\Sigma_{\Delta}$. Define the
{\it extremal critical value of  $e^{-f_0}$ on $\Delta$} as the number
$c(e^{-f_0},\Delta) = e^{-f_0(\Sigma_{\Delta})}$.

Denote  $\Ch(A;f_0)$ the set of all bounded or growing domains of $A$.
Let $\gamma(A) = \left|\Ch(A;f_0)\right|$.  Enumerate these domains by numbers
$1,\ldots,\gamma(A)$.  For every $i\in I$ and $j\in \{1,\ldots ,\gamma
(A) \}$, choose a branch of the multi--valued function $f_i^{\alpha_i}$ on
the domain $\Delta_{j}$ and denote it $g_{i,j}$.  Assume that $|f_i|$ is
bounded on  $\Delta_j$. Let $\Sigma_{i,j}$ be the external support of
$\Delta_j$ with respect to $f_i$. Define the {\it extremal critical
value of the chosen branch $g_{i,j}$ on $\Delta_j$ with respect to
$f_0$} as the number $c(g_{i,j},\Delta_j,f_0)=g_{i,j}(\Sigma_{i,j})$.
 Notice that, if $\Delta_j$ is a bounded domain of $A$, then
$c(g_{i,j},\Delta_j)=c(g_{i,j},\Delta_j,f_0)$.

 Now assume that $|f_i|$ is unbounded on $\Delta_j$. Thus, $\Delta_j$ is a
growing domain of $A$ and $tr(\Delta_j)$ is a bounded face of $tr(A)_{H_0}$.
Denote $M=f_0(\Sigma_{\Delta_j})$. Consider the rational function
$\th=f_i/(f_0-M)$ on $\Delta_j$. Notice that $\th|_{tr(\Delta_j)}$
coincides with the restriction of the linear function $h_i=f_i^0/f_0^0$
to the same set $tr(\Delta_j)$.
Since the sign of $\th$ on $\Delta_j$  is
the same as the sign of $f_i$ on $\Delta_j$ we
can choose a branch of
$\th^{\alpha_i}$ on $\Delta_j$ which has the same argument as $g_{i,j}$
and denote it $\tg$. Let $\Sigma_j$ be the external support of
$tr(\Delta_j)$ with respect to $h_i$ in the affine arrangement
$tr(A)_{H_0}$. Define the {\it the extremal critical
value of the  chosen branch $g_{i,j}$ on $\Delta_j$ with respect to
$f_0$} as the number $c(g_{i,j},\Delta,f_0)=\tg(\Sigma_j)$.

Denote $c(A;\alpha;f_0)$ the product of all extremal critical values  with
respect to $f_0$ of the chosen branches,
 $$ c(A;\alpha;f_0) = \prod_{j=1}^{\gamma(A)}\left(e^{-f_0(\Sigma_{\Delta_j})}
    \prod_{i\in I}c(g_{i,j},\Delta_j,f_0)\right).$$

\section{Hypergeometric period matrix}\label{hg-section}

\subsection{$\beta${\bf nbc}-bases}
\label{bnbcbases}
 Let $A$ be an essential arrangement in an \mbox{$n$--}dimensional real affine
space $V$.
 Define a linear order $<$ in $A$ putting $H_{i}< H_{j}$ if $i<j$.
 A subset $\{H_{i}\}_{i\in J}$ of $A$
is called {\it dependent} if $\cap_{i\in J}H_{i}\neq\emptyset$ and
$\codim(\cap_{i\in J}H_{i})<\left|J\right|$. A
subset of $A$ which has nonempty intersection and is not dependent is
called {\it independent}. Maximal independent sets are called
 {\it bases}. An intersection of a basis defines a vertex.

A $k$-tuple $S=(H_{1},\cdots ,H_{k})$ is called a {\it circuit} if
 $(H_{1},\cdots ,H_{k})$ is dependent and if for each $l$, $1\leq
 l\leq k$, the ($k-1$)-tuple $(H_{1},\cdots ,\widehat{H_{l}},\cdots
 ,H_{k})$ is independent.
A $k$-tuple $S$ is called a {\it broken circuit} if there exists
$H< min (S)$ such that $\{H\}\cup S$ is a  circuit, where $min (S)$
 denotes the minimal element of $S$ for $<$.
\vskip\baselineskip

The collection of subsets of $A$ having nonempty intersection
 and containing no broken circuits is denoted {\bf BC}. {\bf BC}
consists of independent sets. Maximal (with respect to inclusion)
elements of {\bf BC} are bases of $A$ called {\bf nbc}-bases. Recall that
$n$ is the dimension of the affine space.

An {\bf nbc}-basis $B=(H_{i_{1}},\cdots ,H_{i_{n}})$ is called {\it
ordered} if $H_{i_{1}}< H_{i_{2}}< \cdots < H_{i_{n}}$.
The set of all ordered {\bf nbc}-bases of $A$ is denoted  {\bf nbc}$(A)$

A basis $B$ is called a $\bnbc$--{\it{basis}} if $B$ is an
{\bf nbc}-basis and if
\begin{equation} \label{bnbccond}
 \forall H\in B \,\exists H'<\, H \,\mathrm{such}\,\mathrm{ that}\,
  (B-\{H\})\cup \{H'\}\, \mathrm{is}\,\mathrm{ a}\,\mathrm{ base}.
\end{equation}
  Denote $\bnbc (A)$ the set of  all ordered $\bnbc$-bases. Put the
lexicographic order on $\bnbc (A)$
% using the lexicographic order on the hyperplanes read from right to left.

The definition and basic properties of the $\bnbc$-bases are due to
Ziegler [Z].
\vskip\baselineskip

For a basis $B=(H_{i_{1}},\cdots,H_{i_{n}})$, let
$F_{j}=\bigcap_{k=j+1}^{n}H_{i_{k}}$ for $0\leq j\leq n-1$ and
$F_{n} = V$.  Then
$ \xi(B)= (F_{0}\subset F_{1} \subset \cdots\subset F_{n})$ is a flag of
affine subspaces of $V$ with $\dim F_{j} = j$ $(0\leq j\leq n)$.
This flag is called the {\it flag associated with $B$}.

For an edge $F$ of $A$, remind that  $I(F) = \{i\in I \mid F \subseteq
H_{i}\}$. Introduce a differential one-form
$$\omega_{\alpha}(F,A)=\sum_{i\in I(F)} \alpha_{i}\frac{df_i}{f_i}.$$
 For a basis $B=(H_{i_{1}},\cdots,H_{i_{n}})$, let
$\xi(B)= (F_{0}\subset F_{1} \subset \cdots\subset F_{n})$
be the associated flag. Introduce a differential $n$--form $\Xi(B,A)=
\omega_{\alpha}(F_{0},A) \wedge\cdots\wedge\omega_{\alpha}(F_{n-1},A)$.

If $\bnbc (A)=\{B_{1},\cdots,B_{\beta (A)}\}$
and $\phi _{j} = \phi_{j}(A)=\Xi(B_{j},A)$ for $j\in\{1,\ldots,\beta(A)\}$,
define
\begin{equation}
\Phi (A)= \{\phi_{1},\cdots,\phi_{\beta (A)}\}.
\end{equation}
{\sc Example.} For an arrangement $A$ of $p$ hyperplanes in general
position,  the set  $\bnbc(A)$ coincides with the set $\{
(H_{i_1},\ldots,H_{i_n})\, | \, 2\leq  i_1 <\cdots<i_n\leq p\}$. The
latter corresponds to all  vertices of $A$  away from the hyperplane
$H_1$. The differential $n$--forms are
$$ \Phi (A) = \{\alpha_{i_1}\ldots\alpha_{i_n}\frac{df_{i_1}}{f_{i_1}}\wedge
\cdots\wedge\frac{f_{i_n}}{f_{i_n}}\, | \, 2\leq i_1 <\cdots<i_n\leq p\}.$$

\subsection{The definition of the hypergeometric period matrix}
\label{defhg}

Let $\xi = (F_{0} \subset F_{1} \subset \cdots \subset F_{n})$
be a flag of edges of $A$ with $\dim F_{i} = i $ ($i = 0, \ldots,
n-1;$ $F_n=V$).
Let $\Delta$ be a domain of $A$ and $\overline{\Delta} $ its closure in
$V$. We say that the flag is {\it adjacent to} the domain
if $\dim (F_{i} \cap \overline{\Delta}) = i$ for $i = 0,\ldots, n$.

The following proposition from [DT, Proposition 3.1.2] allows us to
enumerate the bounded domains of a configuration $A$ by means of $\bnbc(A)$.

\begin{proposition} \label{bnbclabel}
There exists a unique bijection
$$ C : \bnbc (\A) \longrightarrow \Ch(\A) $$
such that for any $B\in \bnbc(A)$, the associated
flag $\xi(B)$ is  adjacent to the bounded domain $C(B)$.
\end{proposition}

Let $t>0$ be a number which is larger than the maximum of $f_0$ on the
closure of any  bounded domain of $A$. Then the hyperplane $H_t =
\{f_0=t\}$ does not  intersect the bounded domains of $A$. Consider the
affine
arrangement $A_t=A\cup \{H_t\}$. The set of its bounded domains consists of
two disjoint subsets: the first is the  subset of all bounded domains of
$A$; the second is formed by the domains of  $A_t$ which  are
intersections of unbounded domains of $A$ and the  half--space
$\{f_0<t\}$. Notice that the intersection of an unbounded
domain $\Delta$ of $A$ and  the half--space $\{f_0<t\}$ is nonempty and
bounded if and only if $\Delta$ is a growing domain. Thus $\beta (A_t)
=\gamma (A)$.

Define an order $<$ on $A_t$ as $H_t<H_1<\ldots<H_p$. Consider the
set $\bnbc(A_t)$ with respect to this order. If $B \in \bnbc (A_t)$, then
$H_t \not \in B$ because of condition~(\ref{bnbccond}) and the
minimality of $H_t$ with respect to the order $<$.
This observation implies that $\bnbc (A_t)$ and $ \Phi (A_t)$ do not
depend on $t$. Denote them  $\bnbc (A;f_0)$ and  $\Phi (A;f_0)$
respectively. Notice that $\bnbc(A)$ and $\Phi(A)$ are subsets of
$\bnbc(A;f_0)$ and $\Phi(A;f_0)$ respectively because the order on $A$ is
a restriction of the order on $A_t$ to its subset $A$ and because of
condition~(\ref{bnbccond}). We also have an analog of
Proposition~\ref{bnbclabel}.

\begin{proposition} \label{bnbclabel1}
There exists a unique bijection
 $$ \oC : \bnbc (A;f_0) \longrightarrow \Ch(A;f_0) $$
such  that for any $B\in \bnbc (A;f_0)$, the associated flag
 $\xi(B)$ is adjacent to the domain $\oC(B)$.
 Moreover, $\oC|_{\bnbc(A)} = C$. \qed
\end{proposition}

Let the set $\bnbc (A;f_0) = \{B_{1}, \ldots, B_{\gamma}\}$
be lexicographically ordered as in section \ref{bnbcbases}.
For $i = 1, \ldots, \gamma$, define a domain $\Delta_{i} \in \Ch(A;f_0)$
by $\Delta_{i} = C(B_{i})$. This gives us an order on the set of the
growing and bounded domains of $A$. The order is called the
$\bnbc $--{\it{order}}.

We give an orientation to each domain $\Delta\in\Ch(\A;f_0)$ as follows.
 Let $\Delta = C(B)$ with $B\in\bnbc(\A;f_0). $
Let $\xi(B) = (F_{0}\subset F_{1} \subset \cdots\subset F_{n})$
be the associated flag. The flag $\xi(B)$ is adjacent to the domain $B$
and defines its {\it intrinsic orientation} [V2, 6.2]. The
intrinsic orientation is
defined by the unique orthonormal frame $\{e_{1},\ldots, e_{n}\}$ such
that each $e_{i}$ is a unit vector originating  from the point $F_{0}$ in
the direction of $F_{i}\cap \overline{\Delta}$.

Let $\beta = \beta(\A)$.
 Assume that $\Ch(A)=\{\Delta_1, \ldots,\Delta_{\beta}\}$
is the $\bnbc$-ordered set of the bounded domains of $A$  and
$\Phi (\A) = \{\phi_1, \ldots, \phi_{\beta}\}$ is the
$\bnbc$-ordered set of differential $n$--forms constructed  in
Section~\ref{bnbcbases}. Assume that the weights
$\{\alpha_i\}_{i\in I}$ have positive real parts. For every $i\in I$ and
$j\in \{1,\ldots,\beta\}$,
choose a branch of $f_{i}^{\alpha_{i}}$  on the domain $\Delta_{j}$ and
the intrinsic orientation of the domain $\Delta_{j}$.
Let $U_{\alpha} :=f_{1}^{\alpha_{1}} \cdots  f_{p}^{\alpha_{p}}$.
The choice of branches of the functions $f_i^{\alpha_i}$
on all bounded domains defines a choice of branches of the function $U_\alpha$
on all bounded domains.
Define the {\it hypergeometric period matrix} by
\begin{equation} \label{period}
\PM(\A; \alpha) =
\left[\int_{\Delta_j}U_{\alpha}\phi_k\right]_{k,j=1}^{\beta}.
\end{equation}
Since $\Re\alpha_i>0$, all elements of the period matrix are well defined.

Let $\gamma = \gamma(\A)$.
 Let $\Ch(A;f_0)=\{\Delta_1, \ldots,\Delta_{\gamma}\}$
be the $\bnbc$-ordered set of the bounded and growing domains of
$A$  and let $\Phi (\A;f_0) = \{\phi_1, \ldots, \phi_{\gamma}\}$ be the
$\bnbc$-ordered set of differential $n$--forms constructed in
Section~\ref{defhg}. Assume that the weights $\{\alpha_i\}_{i\in I}$ have
positive real parts. For every $i\in I$ and $j\in \{1,\ldots,\gamma\}$,
choose a branch of $f_{i}^{\alpha_{i}}$  on the domain $\Delta_{j}$ and
the intrinsic orientation of each domain $\Delta_{j}$.
The choice of branches of the functions $f_i^{\alpha_i}$
on all bounded and growing
domains defines a choice of branches of the function $U_\alpha$
on all bounded and growing domains.
Define the {\it hypergeometric period matrix with respect to $f_0$} by
\begin{equation} \label{newperiod}
\PM(\A; \alpha; f_0) =
\left[\int_{\Delta_j}e^{-f_0}U_{\alpha}\phi_k\right]_{k,j=1}^{\gamma}.
\end{equation}
Since $\Re\alpha_i>0$ and $f_0$ tends to $+\infty$ on the growing domains
of $A$, all elements of the period matrix are well defined.

\section{The main theorem}\label{mainth}

In [DT], Douai and Terao proved the following theorem, cf. also [V1,V2].

\begin{theorem} \label{dt-main}
Let $A$ be a weighted arrangement given by functions $\{f_i\}_{i\in I}$
and weights $\alpha=\{\alpha_i\}_{i\in I}$ such that $\Re \alpha_{i}>0$ for
all $i\in I$. Fix  branches of the multivalued  functions
$\{f_i^{\alpha_i}\}_{i\in I}$ on all bounded domains of $A$. Then
\begin{equation}\label{dt-formula}
\det \PM (A;\alpha)= c(\A;\alpha) B(A;\alpha).
\end{equation}
\end{theorem}

The main result of this paper is the following theorem.
 \begin{theorem}
\label{maintheorem}
Let $A$ be a weighted arrangement given by functions $\{f_i\}_{i\in I}$
and weights $\alpha=\{\alpha_i\}_{i\in I}$ such that $\Re \alpha_{i}>0$ for
all $i\in I$. Let an additional non-constant linear function $f_0$ be
given. Denote $H_0$ the hyperplane $\{f_0=0\}$. Fix  branches of the
multivalued  functions  $\{f_i^{\alpha_i}\}_{i\in I}$ on all bounded and
growing domains of $A$. Then
\begin{equation}\label{mainformula}
\det \PM (A;\alpha;f_0)= c(\A;\alpha;f_0)B(A;\alpha;H_0).
\end{equation}
\end{theorem}

We will deduce this formula for the determinant of the period matrix with
respect to $f_0$ from Theorem~\ref{dt-main} by passing to a limit.

\section{Proofs} \label{proofs}

\subsection{ Proof of Theorem~\ref{bijection}} \label{biject}

\begin{lemma}
Let $\Delta$ be a growing face. Then
$\overline{tr(\Delta)}\cap \overline{H_0} = \emptyset$, i.e. $tr(\Delta)$
is a bounded face at infinity.
\end{lemma}
$\pf$ Let $x_0,\ldots,x_{n-1}$ be affine coordinates on $V$
such that $f_0(x)=x_0$. Let $(t_0:t_1:\cdots :t_n)$ be the corresponding
projective coordinates in $\overline{V}$: $\{x_i=t_i/t_n\}_{i=0}^{n-1}$.
Let $\Delta$ be a growing face. Assume that $\overline{tr(\Delta)}\cap
\overline{H_0} \not = \emptyset$ and $P=(p_0:p_1:\cdots:p_n)$ is a
point of this intersection. Thus, $p_0=p_n=0$. Let
$Q=(q_0:q_1:\cdots:q_n)$ be any point inside $\Delta$. Thus $q_n\not = 0$.
Since $\overline{\Delta}$ is a closed polyhedron in $\overline{V}$ it
contains the segment $PQ$. This segment is parametrized by the points
$P_{\lambda} = (\lambda p_0 + (1-\lambda)q_0:\lambda p_1 + (1-\lambda)q_1:
\cdots:\lambda p_n + (1-\lambda)q_n)$, for $\lambda \in [0,1]$.
 The point $P_{\lambda}$ tends to $P\in \Hi$  when $\lambda\mapsto 1$. We
have
$$f_0(P_{\lambda})=\frac{\lambda p_0 + (1-\lambda)q_0}{\lambda p_n +
(1-\lambda)q_n}=\frac{q_0}{q_n} = \mathrm{constant}.$$
This contradicts to the assumption that $\Delta$ is a growing face. So
$\overline{tr(\Delta)}\cap \overline{H_0} = \emptyset$. Part $(i)$ of
Theorem~\ref{bijection} is proved. \qed

\begin{lemma}\label{supp-l}
 Let $\Delta$ be an unbounded domain of $A$. Let $tr(\Delta)$ be a bounded
face at infinity with respect to $f_0$.
Let $f_0$ be unbounded on $\Delta\cap\{f_0>0\}$. Then
$\Delta$ is a growing domain of $A$.
\end{lemma}
$\pf$ Since $tr(\Delta)$ is bounded at infinity, we have
$\overline{\Delta}\cap\Hi\cap\overline{H_0}=\emptyset$.
For a real $t$, let $H_t=\{f_0=t\}$. Then
\begin{equation}\label{support-c}
\overline{\Delta}\cap\Hi\cap\overline {H_t}=\emptyset.
\end{equation}
Let $\{x_i\}_{i=1}^{\infty}$ be a sequence of points in $\Delta$ such
that $x_i$ tends to $\infty$ when $i\mapsto\infty$.
 Choose a positive $T$. Assume that $T$ is larger than the supremum of
$f_0$ on all bounded domains of $A$.  Consider the arrangement
$A_T=A\cup\{H_T\}$. Formula~(\ref{support-c}) implies that
$H_T\cap\Delta$ is a bounded domain of the section $(A_T)_{H_T}$. Since
$A_T$ is essential, Proposition~9.9 [BBR] is applicable. It implies that
there is a bounded domain $\Delta_T$ of the arrangement $A_T$, such that
$H_T\cap\Delta$ is a subset of the boundary of $\Delta_T$. This bounded
domain must be $\Delta\cap\{f_0<T\}$, because of the choice of $T$.

Since $\Delta_T$ is bounded, there exists a positive integer $N_T$ such
that for every integer $n\geq N_T$ we have $x_n \in \Delta - \Delta_T$.
Since $\Delta - \Delta_T = \Delta\cap\{f_0\geq T\}$, we have $f_0(x_n)\geq
T$ for all $n\geq N_T$. This proves that $\Delta$ is a growing domain.
\qed

\begin{lemma}  For any bounded face at infinity, $\Sigma$, there exist
exactly $s(\Sigma )$ growing domains with trace $\Sigma$, where $s(\Sigma )$
denotes the discrete width of this face.
\end{lemma}
$\pf$  Let $\Sigma$ be a bounded face at infinity of codimension $k$. Choose
projective coordinates $(t_0:t_1:\cdots:t_n)$ on $\overline{V}$ such that
$\overline{H_0} = \{t_0=0\}$, $\Hi = \{t_n=0\}$, and $F_{\Sigma}$ is
given by $t_1=\cdots =t_{k-1}= t_n= 0$.

Let $v$ be a point in $\Sigma$ and $B$ an open ball around $v$. If the
ball is sufficiently small, then the domains of $A$ which intersect $B$ are
precisely those for which $v$ belongs to their closure in $\overline{V}$
and the
hyperplanes of $\Ai$ which intersect $B$ are exactly those belonging to
$\Ai^{F_{\Sigma}}$. Local affine coordinates on $B$ are given by
$\{y_i=t_i/t_0\}_{i=1}^n$. Since $F_{\Sigma}$ is given by the equations
$y_1=\cdots=y_{k-1}=y_n=0$, the subspace $L$ through $v$ spanned by
the coordinate vectors $e_1,\ldots,e_{k-1},e_n$ is a normal subspace to
$F_{\Sigma}$. Then the number of open domains in $B$ is equal to the
number of open  domains of the arrangement induced in the
tangent space $T_vL$ by the
arrangement
$\Ai^{F_{\Sigma}}$. On $B$ the function $f_0$ has the form $f_0(y)=1/y_n$.
We are interested in  the domains in $B$ on which $f_0 \mapsto +\infty$
when $y_n \mapsto 0$. So, on this domains we must have $y_n >0$. If the
codimension of $\Sigma$ in $\Hi$ is $0$, then the number of such domains
is equal to $1$, which is exactly the discrete width of the empty
configuration. Assume that the above codimension is positive. Then the
number of domains in $B$ on which $y_n>0$ is equal to the number of the
domains of the projective normal arrangement $P\Ai^{F_{\Sigma}}$.
Finally, we want to count only those domains for which $\overline{\Sigma}$
is the only part of their closure in $\overline{V}$, lying in $\Hi$. Thus,
they are the projective domains away from the hyperplane $y_n=0$.
Their number is equal to the discrete length of $P\Ai^{F_{\Sigma}}$. By
definition this number is equal to the discrete width of $\Sigma$.
Lemma~\ref{supp-l} implies that the corresponding domains of $A$ are
growing. \qed

\begin{lemma}
 The number of growing domains of $A$ is equal to the sum of the
volumes of all edges of $\Ai$ at infinity which do not lie in
$\overline{H_0}$.
\end{lemma}
$\pf$ Let $F$ be an edge at infinity with non-zero volume which do not
lie in $\overline{H_0}$. Then, by definition, there are exactly $l(F)$
bounded  faces at infinity which generate $F$. For each of them,
$\Sigma$, there  exist exactly $s(F)$ growing domains of $A$ with trace
$\Sigma$. Thus  there exist exactly $vol(F)=l(F)s(F)$ growing domains
whose traces  generate $F$. Finally, in order to count all growing
domains of $A$,  we have to sum over all edges at  infinity which have
non-zero volume and do not lie in $\overline{H_0}$.
Theorem~\ref{bijection} is proved. \qed

% *********************************************************************

\subsection{ Asymptotic behavior of critical values}\label{critasym}

 Let $A$ be an arrangement in the affine space $V$. Let $f_0$ be an
additional non--constant linear function on $V$.
Define $f_t=1-\frac{f_0}{t}$
% $g_t=(f_t)^t$,
and $H_t=\{f_t=0\}$.
Consider a new weighted
arrangement  $A_t=A\cup \{H_t\}$ where we assume that the weight of
$H_t$ is equal to $t$. For a  sufficiently big $t$, the hyperplane
$H_t$ intersects only some of the unbounded domains of the arrangement
$A$. Moreover, the  intersection creates a new bounded
domain if and only if the intersected  domain  is a growing one. So if
$\Delta$  is a growing domain, we will denote the corresponding bounded
domain of $A_t$ by  $\Delta_t$ and will call it a {\it growing bounded domain}.
If $\Delta$ is a bounded domain of $A$, then it is also a bounded domain of
$A_t$. This correspondence between the bounded domains of $A_t$ and
the bounded or growing domains of $A$ is a bijection.

\begin{lemma}\label{tcrit}
Let $\Delta_t$ be a bounded domain of the arrangement $A_t$. Let $\Delta$
be the corresponding bounded or growing domain of $A$. If $t>0$ and
$(1-f_0/t)$ is positive on $\Delta_t$ choose the positive branch $g_t$ of
$(1-f_0/t)^t$ on $\Delta_t$.

Then the external support of $\Delta_t$ with respect to $f_t$ is a face
of the arrangement $A$. For every big enough $t$ this external support
coincides with the support face, $\Sigma_{\Delta}$, of $f_0$ on $\Delta$.
Moreover, $\lim_{t\mapsto +\infty} c(g_t,\Delta_t) =
e^{-f_0(\Sigma_{\Delta})}$. \end{lemma}

$\pf$ For a fixed $t$, the external support of $\Delta_t$ with  respect
to $f_t$ lies outside $H_t$. Thus, it is a face of the arrangement $A$.

The set of all critical faces of $A$ with respect to $f_0$ is
finite. Let  $M$ be the maximum of $f_0$ on this set. Assume that $t > M$.
Then all bounded domains of $A_t$ lie inside the positive half--space with
respect to $f_t$. Let $f_{t}=1-f_0/t$ attains its maximum
on a critical face $\Sigma_{\Delta}$ of $\overline{\Delta_{t}}$. This
is equivalent to the condition that $f_0$ attains its minimum on the same
face. So $\Sigma_{\Delta}$ is the support face of $\Delta$ with respect
to $f_0$. On the other side, it is the external support of $\Delta_t$
with respect to $f_t$. Hence
$$\lim_{t\mapsto +\infty} c(g_t,\Delta_t) = \lim_{t\mapsto +\infty}
\left(1-\frac{f_0(\Sigma_{\Delta})}{t}\right)^t =
e^{-f_0(\Sigma_{\Delta})}.$$  \qed

\begin{lemma} \label{unbcrit}
  Let the hyperplane $H=\{f=0\}$ belongs to the arrangement $A$.
Let  $\Delta$  be a growing domain of $A$ and  $\Delta_t$ the
corresponding growing bounded domain of $A_t$.  Let $|f|$ be unbounded on
the growing domain $\Delta$.

 Then there exists a unique face  $\Sigma$ of highest dimension,
belonging to  the closure of $\Delta$, such that for every big enough $t$,
the external  support of the face $\Delta_t$  with respect to $f$ is
$\Sigma_t=\Sigma \cap \Delta_t$. Moreover, $tr(\Sigma)$ is the external
support of $tr(\Delta)$ with respect to $h$ in the affine space
$W=\Hi - \overline{H_0}\cap\Hi$, where $h=f^0/f^0_0$. The asymptotic
behavior of $f(\Sigma_t)$ when $t$ tends to $+\infty$ is
given by $f(\Sigma_t)=h(tr(\Sigma))t(1+o(1)).$
\end{lemma}

$\pf$ The set of critical faces with respect to $f$ of the
arrangement $A$ is finite. $|f|$ is bounded on this set. Since $|f|$ is
unbounded on $\Delta$ the external support of $\Delta_t$ with respect to $f$
lies on $H_t$ for $t$ big enough.

 Let $\Sigma_{t_1}$ be a critical face of $\Delta_{t_1}$ which lies on
$H_{t_1}$ for some $t_1$ fixed. Then $\Sigma_{t_1} = H_{t_1}\cap\Sigma$,
where $\Sigma$ is a face of $\overline{\Delta}$. Consider the face
$\Sigma_t=H_t\cap\Sigma$ of $\Delta_t$ for an arbitrary $t$. It is a
critical face of $\Delta_t$ because $\Sigma_t$ is parallel to
$\Sigma_{t_1}$ and the latter is parallel to $H$. Let us compute the
asymptotic behavior of $f(\Sigma_t)$ when $t$ tends to $+\infty$.

Choose affine coordinates $\{x_j\}_{j=0}^{n-1}$ in $V$ such that
$f_0(x)=x_0$. Let $f=f^0 + b$ be the sum of the homogeneous part of $f$ and
the constant term. Then $f=x_0(f^0/x_0 + b/x_0) = x_0(f^0/f_0^0 + b/x_0)$.
Since $H_t=\{x_0=t\}$,  $f(\Sigma_t) = t(h(\Sigma_t) + b/t) =
h(tr(\Sigma))t(1+o(1))$.

Let $\Sigma'$ be the external support of $tr(\Delta)$ relative to $tr(H)$.
Let $\Sigma$ is the face of $\Delta$ for which $\Sigma'=tr(\Sigma)$. Then
the previous computation shows that for every $t$ big enough $\Sigma_t$
is the external support of $\Delta_t$ with respect to  $f$. \qed

\begin{corollary}\label{unb-cor}
Let the conditions be as in Lemma~\ref{unbcrit}. In addition, assume
that $\alpha$ is a complex number. Fix a branch of $f^{\alpha}$ on
$\Delta$ and denote it $g$. Fix a branch of $(f/f_0)^{\alpha}$ on $\Delta$
as in Section~\ref{critical} and denote it $\tilde{g}$. Fix branches
of $t^{\alpha}$ and $(1+o(1))^{\alpha}$ using the branch of the
logarithm with zero argument. Then the
asymptotic behavior of $c(g,\Delta_t)$ when $t$ tends to $+\infty$ is
$c(g,\Delta_t)=c(g,\Delta,f_0)t^{\alpha}(1+o(1))$.
\end{corollary}
$\pf$ Use the notation of the previous proof. Since
$c(g,\Delta_t)=g(\Sigma_t)$, $c(g,\Delta,f_0)=\tilde{g}(tr(\Sigma))$ and
the arguments of $g$ and $\tilde{g}$ are the same on $\Delta$, the
assymptotic formula for $f$ implies the statement of the corollary.\qed

\begin{lemma} \label{bndcrit}
  Let the hyperplane $H=\{f=0\}$ belongs to the arrangement $A$. Let
$\alpha$ be a complex number. Let  $\Delta$  be a  growing domain of $A$
and  $\Delta_t$ the  corresponding growing bounded domain of $A_t$.  Fix a
branch of $f^{\alpha}$ on $\Delta$ and denote it $g$. Then

 $|f|$ is bounded on $\Delta$ if and only if $tr(\Delta)\subset
tr(H)$. The latter condition is equivalent to the equation
$h(tr(\Delta))=0$ where  $h=f^0/f^0_0$ is a linear function of
the arrangement $tr(A)_{H_0}$. Moreover, if $|f|$ is bounded
on $\Delta$, then for every big enough $t$, $c(g,\Delta_t)$ equals
 the constant $c(g,\Delta,f_0)$.
\end{lemma}

$\pf$ $|f|$ is bounded on $\Delta$ if and only if $\Delta$ is placed between
two hyperplanes $H'$ and $H''$ parallel to $H$. Denote the domain between
these two hyperplanes by $D$. Since $\Delta\subset D$ we have
$tr(\Delta)\subset tr(D)=tr(H)$.

The reverse part is a consequence of Lemma~\ref{unbcrit}.

 Let $\Sigma$ be the external support of $\Delta$ with respect to  $f$.
Then for every big enough $t$, $c(g,\Delta_t)=g(\Sigma)$.
The latter equals $c(g,\Delta,f_0)$. \qed

% ********************************************************************

\subsection{ Proof of Theorem~\ref{maintheorem}}

We prove Theorem~\ref{maintheorem}  applying Theorem~\ref{dt-main}  to
the arrangement $A_t$ and then passing to the limit when $t\mapsto +\infty$.

First study $B(A_t;\alpha,t)$.
\begin{lemma} \label{b-asym}
$(i)$ The only factor in the numerator of $B(A_t;\alpha,t)$ depending on
$t$, when $t$ is big enough, is the factor corresponding to the edge $H_t$.
It contributes $\Gamma(t+1)^{\mathrm{\#}},$ where \# is the number
 of growing domains of $A$.

$(ii)$ The factors in the denominator depending on $t$ come from the
edges at infinity with non-zero volume which do not lie in $\overline{H_0}$.
Each of them, $F$, contributes $\Gamma(t+1+\alpha'(F))^{vol(F)}$, where
$\alpha'(F) = \sum_{H\in A;\, F\not\subseteq \overline{H}} \alpha_H$.

$(iii)$ The asymptotic behavior of $B(A_t;\alpha,t)$ when $t$ tends to
$+\infty$ is given by
$$ B(A_t;\alpha,t)=B(A;\alpha;H_0)\prod_{F\in L_{-};\, F\not\subset
\overline{H_0}}t^{-\alpha'(F)vol(F)}(1+o(1)).$$
\end{lemma}

$\pf$ Recall that
$$B(A_t;\alpha,t) = \prod_{F\in L_{t+}} \Gamma (\alpha (F) +1)^{vol(F)}
 \big/   \prod_{F\in L_{t-}} \Gamma (-\alpha (F) +1)^{vol(F)} , $$
$$B(A;\alpha ;H_0) = \prod_{F\in L_{+}} \Gamma (\alpha (F) +1)^{vol(F)}
 \big/  \prod_{F\in L_{-}; F\subset H_0} \Gamma (-\alpha (F) +1)^{vol(F)} ,$$
 where
$L_{t-}$, $L_{-}$ denote the set of all edges at infinity of the
arrangements $\overline{A_t}$ and $\Ai$, respectively, and
$L_{t+}$, $L_{-}$ denote the set consisting of all the other edges of
the same arrangements.

 $(i)$  Since the only weight depending on $t$
corresponds to $H_t$, the factors in the denominator that depend on $t$
correspond to the edges of $A_t$ lying in $H_t$. If such an edge $F$ is a
proper subspace of $H_t$, then it is decomposable [STV, Section 2], that is
the localization of the arrangement $A_t$ at the edge $F$ is a product of
two nonempty subarrangements where one of the subarrangements is equal to
$\{H_t\}$. According to [STV, Proposition 7], the discrete width of a
decomposable edge is zero. Thus its discrete volume is zero.

 The volume  of $H_t$ is the number of bounded domains of the section
arrangement  $(A_t)_{H_t}$ which is exactly the number of growing domains
of the arrangement $A$.

$(ii)$ If $F$ is an edge at infinity of $A_t$, then $\alpha(F) =
-t-\sum_{i\in I}\alpha_p + \sum_{H\in A_t;\, F\subset
\overline{H}}\alpha_H$.
The last  sum depends on $t$ if and only if $\sum_{H\in A_t;\, F\subset
\overline{H}}\alpha_H$  does not depend on $t$, i.e. if and only if
$F\not\subset \overline{H_t}$. Since  $\overline{H_t}\cap H_{\infty} =
\overline{H_0}\cap H_{\infty}$, the weight  $\alpha(F)$ depends on
$t$ if and only if $F\not\subset \overline{H_0}$. So  $\alpha(F)= -t
-\alpha'(F)$.
Notice that such an edge is also an edge of the arrangement $A$.

$(iii)$ According to Theorem~\ref{bijection} the number of growing domains
of the arrangement $A$ is equal to the sum of the volumes of all edges at
infinity of the arrangement $A$ which do not lie in $\overline{H_0}$. Thus
the number of factors in the numerator and in the denominator containing
$t$ is equal. Sterling's formula gives us
$\Gamma(t+1)/\Gamma(t+1+a)=t^{-a}(1+o(1))$ when $t$ tends to $+\infty$.
So we obtain the required formula. \qed
\vskip\baselineskip

  Now consider the limit of the product of the critical values,
$c(A_t;\alpha,t)$. For every $i\in I$ and every bounded or growing domain
$\Delta$ of the arrangement $A$, choose a branch of $f_i^{\alpha_i}$ on
$\Delta$ and denote it $g_{i,\Delta}$. This also fixes
branches of $f_i^{\alpha_i}$ on the bounded domains of $A_t$ independently
on $t$. Notice that for every big enough $t$, $\,f_t$ is positive on all
bounded domains of the arrangement $A_t$. Choose the positive branch of
$(f_t)^t$ on this domains and denote it $g_t$.

\begin{lemma} \label{c-asym}
$c(A_t;\alpha,t)$ has the following asymptotic behavior when $t$
tends to $+\infty$:
 $$ c(A_t;\alpha,t)=c(A;\alpha;f_0)\prod_{F\in L_{-}; F\not\subset
H_0}t^{\alpha'(F)vol(F)} (1+o(1)).$$
\end{lemma}

$\pf$
\begin{eqnarray}
 C(A_t;\alpha;t) & = & \prod_{\Delta\in Ch(A_t)}\left( c(g_t,\Delta)
  \prod_{i\in I} c(g_{i,\Delta},\Delta)\right) \nonumber\\
  & = & \left(\prod_{\Delta\in Ch(A_t)}c(g_t,\Delta)\right)
\left(\prod_{\Delta\in Ch(A)}\prod_{i\in I}c(g_{i,\Delta},\Delta)\right)
\left(\prod_{\Delta_t}\prod_{i\in I}
      c(g_{i,\Delta_t},\Delta_t)\right),\label{3prod}
\end{eqnarray}
where $\Delta_t$ in the last product ranges over the growing bounded
domains of $A_t$.

Describe the asymptotic behavior of each of the three products in
formula~(\ref{3prod}).

Assume that $\Delta$ is a bounded domain of $A_t$. Lemma~\ref{tcrit} asserts
that $\lim_{t\rightarrow +\infty}c(g_t,\Delta)=e^{-f_0(\Sigma_{\Delta'})}$,
where $\Delta'$ is the domain of $A$ (bounded or growing) which corresponds
to the domain $\Delta$ of $A_t$ and $\Sigma_{\Delta'}$ is the
support face of $f_0$ on  $\Delta'$.

If $\Delta$ is a bounded domain of $A$ and $i\in I$, then
$c(g_{i,\Delta},\Delta)=c(g_{i,\Delta},\Delta,f_0)$ by definition.

Let $\Delta_t$ be a growing bounded domain of $A_t$ and $\Delta$ the
corresponding growing domain of $A$. Let $i\in I$. If
$tr(\Delta)\not\subset \overline{H_i}$, then
$c(g_{i,\Delta_t},\Delta_t)=c(g_{i,\Delta},\Delta,f_0)t^{\alpha_i}(1+o(1))$,
by Corollary~\ref{unb-cor}. If $tr(\Delta)\subset \overline{H_i}$, then
$c(g_{i,\Delta_t},\Delta_t)=c(g_{i,\Delta},\Delta,f_0)$ by
Lemma~\ref{bndcrit}.
Let $F$ be an edge at infinity. Assume that $F$ does not lie in
$\overline{H_0}$ and has a non-zero volume. According to
Theorem~\ref{bijection}, there exist exactly $vol(F)$ growing domains of
the arrangement $A$, whose traces generate $F$. Every term in the last
product of formula~(\ref{3prod}) depends on a growing bounded domain.
Collect all the terms such that the trace of the corresponding growing
domain generates $F$. Then for the product of the chosen factors we have
\begin{eqnarray*}
\prod_{\dind}\prod_{i\in I} c(g_{i,\Delta_t},\Delta_t) &=&
  \left( \prod_{\dind}\prod_{i\in I}c(g_{i,\Delta},\Delta,f_0)\right)
   \left( \prod_{\dind}\prod_{i,\, F\not\subset H_i}
     t^{\alpha_i}(1+o(1))\right)\\
 &=&  \left(\prod_{\dind}\prod_{i\in I} c(g_{i,\Delta},\Delta,f_0)\right)
    \left(\prod_{\dind} t^{\alpha'(F)}(1+o(1))\right)\\
 &=&       t^{\alpha'(F)vol(F)}(1+o(1))
      \prod_{\dind}\prod_{i\in I} c(g_{i,\Delta},\Delta,f_0)
\end{eqnarray*}

Collecting the  asymptotic behavior for the three products in
(\ref{3prod}), we obtain the statement of the lemma.\qed
\vskip\baselineskip

{\sf Proof of Theorem~\ref{maintheorem}:} Apply
formula~(\ref{dt-formula}) to the weighted arrangement $A_t$.
Lemmas~\ref{b-asym} and \ref{c-asym}  show that the terms dependent on
$t$ in the asymptotic formulas for $B(A_t;\alpha,t)$ and
$c(A_t;\alpha,t)$ cancel out. Thus,
$$\lim_{t\mapsto +\infty} c(A_t;\alpha,t) B(A_t;\alpha,t)=
   c(A;\alpha;f_0) B(A;\alpha;H_0),$$
which gives us the right hand side of formula~(\ref{mainformula}).

  Let us study the entries of the period matrix $\PM (A_t;\alpha,t)$:
$\,\PM_{k,j}(t)=\int_{\Delta_j}U_{\alpha,t}\phi_k(A_t)$,
where $\Delta_j$ is a bounded domain of $A_t$ and $\phi_k(A_t)$ is one of
the $n$--forms constructed in Section~\ref{bnbcbases}. Remind that for a
fixed $k$ and a
big enough $t$,
the form $\phi_k(A_t)$ is independent on $t$ and equals
$\phi_k(A;f_0)$.  Since $U_{\alpha,t}=(1-f_0/t)^tU_{\alpha}$, we have
$\lim_{t\mapsto +\infty} U_{\alpha,t} = e^{-f_0}U_{\alpha}$.

 Since $\Delta_j$ is a bounded domain of $A_t$, there exists a unique
bounded or growing domain, $\Delta$, of $A$ such that $\Delta_j=\Delta\cap
\{f_0<t\}$. Extend $f_t$ as zero on $\Delta - \Delta_j$. Then
$\PM_{k,j}(t)=\int_{\Delta}(f_t)^tU_{\alpha}\,\phi_k(A;f_0)$. Since
$(f_t)^t<e^{-f_0}$ on $\Delta\cap\{f_0>0\}$, Lebesgue's convergence
theorem is applicable and  $\lim_{t\mapsto +\infty}\PM_{k,j}(t) =
\int_{\Delta}e^{-f_0}U_{\alpha}\,\phi_k(A;f_0)= \PM_{k,j}(A;\alpha;f_0)$.
These limits give us  $\lim_{t\mapsto +\infty}\PM (A_t;\alpha,t)
=\PM (A;\alpha;f_0)$. Theorem~\ref{maintheorem} is proved.\qed

\section{Determinant formulas for Selberg type integrals}\label{selberg}

Let $z_1<\cdots<z_p$ be real numbers. Let
$\alpha_1,\ldots,\alpha_p,\gamma$ be complex numbers with positive real
parts. For $t\in \mathbb{R}^n$ define
$$
\Phi(t,z) = \prod_{s=1}^p\prod_{i=1}^n(t_i-z_s)^{\alpha_s}\prod_{1\leq
i<j\leq n} (t_j-t_i)^{2\gamma}.$$
The branches of $x^{\alpha_s}$ and $x^{2\gamma}$
are fixed by $-\pi/2< \arg\,x < 3\pi/2$ for all $s\in\{1,\ldots,p\}$.

 Let $\mathcal{Z}_n^p=\{ \l=(l_1,\ldots,l_p)\in \mathbb{Z}^p\,|\, l_i\geq
0,\,\,\,\, l_1+\cdots+l_p=n\}$. For every $s\in\{1,\ldots,p\}$ denote
$\l^s=\sum_{i=1}^{s} l_i$, $\l^0=0$. Let $\m \in \mathcal{Z}_n^p$ and $s\in
\{1,\ldots,p\}$. Denote $\Gamma_{\m,s}$ the set of integers
$\{\m^{s-1}+1,\ldots,\m^s\}$ and $d^nt=dt_1\wedge\ldots\wedge dt_n$.
Define the following $n$--forms
$$\omega_{\m}(t,z) =\left(
\sum_{\sigma\in\mathbb{S}^n}\prod_{s=1}^{p}\frac{1}{m_s!}
\prod_{j\in\Gamma_{\m,s}}\frac{1}{(t_{\sigma_j}-z_s)}\right) d^nt.$$
If $\m\in \mathcal{Z}_n^{p-1}$, then we identify $\m$ with the $p$--tuple
$(\m,0)\in \mathcal{Z}_n^p$.

 For $\l\in\mathcal{Z}_n^{p-1}$, let
$$\mathbb{U_{\l}}=\{t=(t_1,\ldots,t_n)\in\mathbb{R}^n
\,|\,z_s\leq t_{\l^{s-1}+1}\leq\cdots\leq t_{\l^s}\leq z_{s+1} \mbox{ for
all } s=1,\ldots,p-1\}.$$
Assume that all domains in the formulas below inherit the standard
orientation from $\mathbb{R}^n$.

\begin{theorem}\label{no-exp}
$\mathrm{,}\,\,\mathrm{cf}\,\,\mathrm{[V6].}$
\begin{eqnarray}
%%!!
\lefteqn{\det\!\left[
\int_{\mathbb{U}_{\l}}\Phi(t,z)\omega_{\m}(t,z)
\right]_{\l,\m\in\mathcal{Z}_n^{p-1}} = } \label{no-exp-form}\\
 & &\prod_{s=0}^{n-1}\left[
\frac{\Gamma((s+1)\gamma)^{p-1}}{\Gamma(\gamma)^{p-1}}
\frac{\Gamma(1+\alpha_p+s\gamma)\prod_{j=1}^{p-1}\Gamma(\alpha_j+s\gamma)}
     {\Gamma(1+\sum_{j=1}^p\alpha_j +(2n-2-s)\gamma)}
\right]^{\binom{p+n-s-3}{p-2}}\nonumber\\
 & & \exp\Bigl(i\pi\binom{p+n-2}{p-1}\sum_{s=1}^{p} (s-1)\alpha_s\Bigr)
 \prod_{1\leq a<b\leq p} (z_b-z_a)^{(\alpha_a+\alpha_b)
 \binom{p+n-2}{p-1}+2\gamma\binom{p+n-2}{p}}. \nonumber
\end{eqnarray}
\end{theorem}

Notice, that formula~(\ref{no-exp-form}) is not symmetric with respect to
$\alpha_1,\ldots,\alpha_p$. To make it symmetric we  introduce new
differential $n$--forms, $\widetilde{\omega_{\m}}$, for
$\m\in\mathcal{Z}_n^{p-1}$. Namely
\begin{eqnarray*}
\widetilde{\omega_{\m}}(t,z) & = &
[\prod_{s=1}^{p-1}(m_s!)\alpha_s(\alpha_s+\gamma)
\ldots(\alpha_s+(m_s-1)\gamma)]\omega_{\m}\\
 & = & \sum_{\sigma\in\mathbb{S}^n}\prod_{s=1}^{p-1}
\alpha_s(\alpha_s+\gamma)\ldots(\alpha_s+(m_s-1)\gamma)
\prod_{j\in\Gamma_s} \frac{1}{(t_{\sigma_j}-z_s)}d^{n}t.
\end{eqnarray*}

Theorem~\ref{no-exp} implies
\begin{eqnarray}
%%!!
\lefteqn{\det\!\left[
\int_{\mathbb{U}_{\l}}\Phi(t,z)\widetilde{\omega_{\m}}(t,z)
\right]_{\l,\m\in\mathcal{Z}_n^{p-1}} = }\label{no-exp-form1}\\
 & &\prod_{s=0}^{n-1}\left[
\frac{\Gamma((s+1)\gamma+1)^{p-1}}{\Gamma(\gamma+1)^{p-1}}
\frac{\prod_{j=1}^{p}\Gamma(\alpha_j+s\gamma+1)}
     {\Gamma(1+\sum_{j=1}^p\alpha_j +(2n-2-s)\gamma)}
\right]^{\binom{p+n-s-3}{p-2}} \nonumber \\
 & & \exp\Bigl(i\pi\binom{p+n-2}{p-1}\sum_{s=1}^{p} (s-1)\alpha_s\Bigr)
 \prod_{1\leq a<b\leq p}
(z_b-z_a)^{(\alpha_a+\alpha_b)\binom{p+n-2}{p-1}+2\gamma\binom{p+n-2}{p}}.
\nonumber
\end{eqnarray}

\begin{lemma}\label{eq-crit}
  For every $\l\in\mathcal{Z}_n^{p-1}$, $i,j \in \{1,\ldots,n\}$,
$s\in\{1,\ldots.p\}$,
fix  branches $g_{j,s}$, $h_{j,i}$ of
the multivalued functions $(t_j-z_s)^{\alpha_s}$ and $(t_j-t_i)^{2\gamma}$
,respectively, on the domain $\mathbb{U}_{\l}$ as at the beginning of
the current section. Then the product
 $$ \exp\Bigl(i\pi\binom{p+n-2}{p-1}\sum_{s=1}^{p} (s-1)\alpha_s\Bigr)
\prod_{1\leq a<b\leq p}
(z_b-z_a)^{(\alpha_a+\alpha_b)\binom{p+n-2}{p-1}+2\gamma\binom{p+n-2}{p}}$$
equals the product of critical values of the chosen branches
\begin{equation}\label{prod-no-exp}
 \prod_{\l\in\mathcal{Z}_n^{p-1}}\left[ \prod_{j=1}^n\prod_{s=1}^p
c(g_{j,s},\mathbb{U}_{\l})\prod_{1\leq i<j \leq n}
c(h_{j,i},\mathbb{U}_{\l})\right].
\end{equation}
The critical values were defined in Section~\ref{critical}.
\end{lemma}

  Lemma~\ref{eq-crit} allow us to replace the last lines in
formulas~(\ref{no-exp-form}) and (\ref{no-exp-form1}) by the product of
critical values~(\ref{prod-no-exp}), cf. [V6].

 For $\l \in \mathcal{Z}_n^p$, let $z_0=-\infty$ and
$$\widetilde{\mathbb{U_{\l}}}=\{t=(t_1,\ldots,t_n)\in\mathbb{R}^n
\,|\,z_{s-1}\leq t_{\l^{s-1}+1}\leq\cdots\leq t_{\l^s}\leq z_s \mbox{ for
all } s=1,\ldots,p\}.$$

\begin{theorem}\label{expon}
 Let $a$ be a complex number with positive real part. Then
\begin{eqnarray}
%%!!
\lefteqn{\det\!\left[
\int_{\widetilde{\mathbb{U}_{\l}}}\exp\Bigl(a\sum_{j=1}^n t_j\Bigr)
\Phi(t,z)\omega_{\m}(t,z)
\right]_{\l,\m\in\mathcal{Z}_n^p} = (-1)^{n\binom{p+n-1}{p-1}}
 }\label{expon-form}\\
 & &
 \prod_{s=0}^{n-1}\left[\frac{\Gamma((s+1)\gamma)^p}{\Gamma(\gamma)^p}
\prod_{j=1}^p\Gamma(\alpha_j+s\gamma)\right]^{\binom{p+n-s-2}{p-1}}
 \prod_{1\leq a<b\leq p}
(z_b-z_a)^{(\alpha_a+\alpha_b)\binom{p+n-1}{p}+2\gamma\binom{p+n-1}{p+1}}
\nonumber \\
 & &  \exp\Bigl(i\pi\Bigl[\binom{p+n-1}{p}\sum_{s=1}^{p} s\alpha_s\Bigr]\Bigr)
 \exp\Bigl(a\pi\binom{p+n-1}{p}\sum_{s=1}^{p} z_s\Bigr)\,
 a^{-\binom{p+n-1}{p}\sum_{s=1}^{p}s\alpha_s-2p\binom{p+n-1}{p+1}\gamma}
\nonumber
\end{eqnarray}
\end{theorem}

  The next lemma allow us to replace the last line in
formula~(\ref{expon-form}) by the product of
critical values~(\ref{prod-expon}).

\begin{lemma}\label{eq-crit1}
  For every $\l\in\mathcal{Z}_n^p$, $i,j \in \{1,\ldots,n\}$,
$s\in\{1,\ldots.p\}$,
fix  branches $g_{j,s}$, $h_{j,i}$ of
the multivalued functions $(t_j-z_s)^{\alpha_s}$ and $(t_j-t_i)^{2\gamma}$
respectively, on the domain $\widetilde{\mathbb{U}_{\l}}$ as in the
beginning of the current section. Then the product
 $$\exp\Bigl(i\pi[\binom{p+n-1}{p}\sum_{s=1}^{p} s\alpha_s]\Bigr)
 \exp\Bigl(a\pi\binom{p+n-1}{p}\sum_{s=1}^{p} z_s\Bigr)\,
 a^{-\binom{p+n-1}{p}\sum_{s=1}^{p}s\alpha_s-2p\binom{p+n-1}{p+1}\gamma}$$
equals the product of critical values of the chosen branches with respect
to the linear function $-at_1$.
\begin{equation}\label{prod-expon}
\prod_{\l\in\mathcal{Z}_n^p}\left[
c(e^{a\sum_{j=1}^nt_j},\widetilde{\mathbb{U}_{\l}})
\prod_{j=1}^n\prod_{s=1}^p
c(g_{j,s},\widetilde{\mathbb{U}_{\l}},-at_1)\prod_{1\leq i<j \leq n}
c(h_{j,i},\widetilde{\mathbb{U}_{\l}},-at_1)\right].
\end{equation}
The critical values were defined in Section~\ref{critical}.
\end{lemma}

\section{Proofs of Theorem~\ref{no-exp} and Theorem~\ref{expon}}
\label{pfselb}

 Theorem~\ref{no-exp} is a direct corollary of Theorems~5.15 and 7.8 [TV].
The computations are long but straightforward.\qed

  The correspondence in notation between the current paper and [TV] is as
follows:
\newpage
\begin{table}[h]
\begin{tabular}{l|c|c}
Object & current notation & [TV] -- article\\
\hline
Dimension of the vector space & n & l\\
Number of points & p & n\\
Coordinates      & $t$  & $u$ \\
Weights          & $\alpha$ & $2\Lambda/p$ \\
                 & $\gamma$   & $-1/p$ \\
Points           & $z \in \mathbb{R}$ & $y/h=z \in i\mathbb{R}$\\
Parameter        & a          & $i\eta/p$
\end{tabular}
\end{table}

\vskip\baselineskip

 For $\l\in\mathcal{Z}_n^p$, let $z_0=-\infty$ and
$$\mathbb{V_{\l}}=\{t=(t_1,\ldots,t_n)\in\mathbb{R}^n
\,|\, z_{s-1}\leq t_{j_s}\leq z_{s} \mbox{ for all } s=1,\ldots,p \mbox{
and } j_s\in\Gamma_{\l,s} \}.$$
Theorems~5.15 and 7.6 [TV] imply the following formula
\begin{eqnarray}
%%!!
\lefteqn{\det\!\left[
\int_{\mathbb{V}_{\l}}\exp\Bigl(a\sum_{j=1}^n t_j\Bigr)
\prod_{s=1}^p\prod_{i=1}^n(t_i-z_s)^{\alpha_s}\prod_{1\leq
i<j\leq n} (t_i-t_j)^{2\gamma}\omega_{\m}(t,z)
\right]_{\l,\m\in\mathcal{Z}_n^p} = }\label{rect-exp} \\
 & & (-1)^{n\binom{p+n-1}{p-1}}
\prod_{\l\in\mathcal{Z}_n^p}\prod_{j=1}^p\prod_{s=1}^{l_j}
\frac{\sin(-s\pi\gamma)}{\sin(-\pi\gamma)}
\prod_{s=0}^{n-1} \left[\frac{\Gamma((s+1)\gamma)^p}{\Gamma(\gamma)^p}
\prod_{j=1}^p\Gamma(\alpha_j+s\gamma)\right]^{\binom{p+n-s-2}{p-1}}
\nonumber \\
& & \exp\Bigl(a\pi\binom{p+n-1}{p}\sum_{s=1}^{p} z_s\Bigr)\,
a^{-\binom{p+n-1}{p}\sum_{s=1}^{p}s\alpha_s-2p\binom{p+n-1}{p+1}\gamma}
\nonumber \\
 & & \exp\Bigl(i\pi\Bigr[\binom{p+n-1}{p}\sum_{s=1}^{p} s\alpha_s
+ p^2\binom{p+n-1}{p+1}\gamma\Bigr]\Bigr)
 \prod_{1\leq a<b\leq p}
(z_b-z_a)^{(\alpha_a+\alpha_b)\binom{p+n-1}{p}+2\gamma\binom{p+n-1}{p+1}}.
\nonumber
\end{eqnarray}

In order to obtain Theorem~\ref{expon} we have to pass from
"rectangular" domains $\mathbb{V}_{\l}$ to "triangular" domains
$\widetilde{\mathbb{U}_{\l}}$. For any
$\l,\,\m \in\mathcal{Z}_n^p$ we have
\begin{eqnarray*}
\lefteqn{\int_{\mathbb{V}_{\l}}\exp\Bigl(a\sum_{j=1}^n t_j\Bigr)\prod_{s=1}^p
\prod_{i=1}^n(t_i-z_s)^{\alpha_s}\prod_{1\leq i<j\leq n}
(t_i-t_j)^{2\gamma}\omega_{\m}= } \\
& = &
\sum_{\sigma\in\mathbb{S}^{l_1}\times\cdots\times\mathbb{S}^{l_p}}
\int_{\sigma\widetilde{\mathbb{U}_{\l}}}
\exp\Bigl(a\sum_{j=1}^n t_j\Bigr)\prod_{s=1}^p
\prod_{i=1}^n(t_i-z_s)^{\alpha_s}\prod_{1\leq i<j\leq n}
(t_i-t_j)^{2\gamma}\omega_{\m} \\
& = & e^{n(n-1)i\pi\gamma}
[\prod_{j=1}^p1(1+e^{-2\pi i\gamma})\cdots(1+e^{-2\pi i\gamma} +\cdots +
e^{-2\pi i \gamma(l_j-1)})]\\
& &
\int_{\widetilde{\mathbb{U}_{\l}}}
  \exp\Bigl(a\sum_{j=1}^nt_j\Bigr)\Phi(t,z)\omega_{\m}(t,z)\\
&=& e^{n(n-1)i\pi\gamma}
\prod_{j=1}^p\prod_{s=1}^{l_j}\frac{\sin(-s\pi\gamma)}{\sin(-\pi\gamma)}
\prod_{j=1}^p e^{-i\pi\gamma l_j(l_j-1)/2}
\int_{\widetilde{\mathbb{U}_{\l}}}
\exp\Bigl(a\sum_{j=1}^nt_j\Bigr)\Phi(t,z)\omega_{\m}(t,z)
\end{eqnarray*}
This proves Theorem~\ref{expon}. \qed


\begin{thebibliography}{999}

\bibitem[BBR]{BBR}
Barnabei, M., Brini, A., Rota, Dgh.-K.[G.-C]: The theory of M\"{o}bius
functions, Uspekhi Mat. Nauk. {\bf 41}, no.3, 113-157 (1986); English
transl. in Russian Math. Surveys {\bf 41} (1986).

\bibitem[DT]{DT}
Douai, A., Terao, H.:  The determinant of a hypergeometric period matrix,
Invent. math. {\bf 128}, 417-436 (1997).

\bibitem[FT]{FT}
Falk, M. J., Terao, H.:  $\bnbc$--bases for  cohomology of local system
on hyperplanes complements, Transactions AMS, {\bf 349}, 189-202 (1997)

\bibitem[L]{L} Loeser, F.: Arrangements d'hyperplans et somme de Gauss,
Ann. Sci. \'Ecole Norm. Sup.
{\bf 24}, 379-400 (1991)

\bibitem[LS]{LS} Loeser, F. and Sabbah, C.:
 Equations aux diff\'erences finies et
d\'eterminants d'int\'egrales de fonctions multiformes, Comment. Math.
Helv. {\bf 66}, 458-503 (1991)

\bibitem[OT]{OT1} Orlik, P., Terao, H.: Arrangements of hyperplanes.
Grundlehren der Math. Wiss. {\bf 300}, Berlin Heildelberg New York:
Springer, 1992.

\bibitem[STV]{STV} Schechtman, V., Terao, H., Varchenko, A.:
Local systems over complements of hyperplanes and the Kac-Kazhdan
conditions for singular vectors, J. Pure Appl. Algebra {\bf 100}, 93-102
(1995).

\bibitem[SV]{ScV1} Schechtman, V., Varchenko, A.: Arrangements of
hyperplanes and Lie algebra homology, Invent. math. {\bf 106},
139--194 (1991).

\bibitem[TV]{TV} Tarasov, V., Varchenko, A.:
Geometry of $q$--hypergeometric functions as a bridge between Yangians and
quantum affine algebras, Invent. Math {\bf 128}, 501-588 (1997).

\bibitem[V1,V2]{Var}
Varchenko, A.: The Euler Beta-function, the Vandermonde determinant,
Legendre's equation, and critical values of linear functions on a configuration
of hyperplanes, Math. USSR Izvestija {\bf 35} (1990), 543-572 and {\bf 36},
 155-168 (1991).

\bibitem[V3]{Var3} Varchenko, A.: Multidimensional hypergeometric
functions and  representation theory of  Lie algebras
and quantum groups, Advanced Series in Mathematical Physics -
{\bf 21},  World Scientific Publishers, Singapore, 1995.

\bibitem[V4]{Var4}
Varchenko, A.: Critical values and the deteminant of the periods,
Uspekhi Mat. Nauk. {\bf 44}, no.4, 235-236 (1989); English
transl. in Russian Math. Surveys {\bf 44} no.4, 209-210 (1989).

\bibitem[V5]{Var5}
Varchenko, A.: Critical points of the product of powers of linear
functions and families of bases of singular vectors,
Compositio Math.  {\bf 97}, 385-401 (1995).

\bibitem[V6]{Var6}
Varchenko, A.: Determinant formula for Selberg-type integrals,
Func.Anal.Appl.  {\bf 25}, no. 4, 304-305 (1991).

\bibitem[Z]{Z} Ziegler, G.: Matroid shellability, $\beta$--systems,
and affine arrangements, J. Alg. Combinatorics, {\bf 1},
283--300 (1992).

\end{thebibliography}
\end{document}